\documentclass[12pt]{article}
\makeatletter
\def\@begintheorem#1#2{\it \trivlist
\item[\hskip \labelsep{\sc #1\ #2.}]}
\makeatother
\usepackage{latexsym}
\usepackage{amsfonts}
\usepackage{amsmath}
\usepackage{amsthm}
\usepackage{amssymb}
\usepackage{amscd}

\newcommand{\us}{u_{\sigma}}

\newcommand{\ut}{u_{\tau}}

\newcommand{\up}{u_{\pi}}
\newcommand{\Up}{u_{\pi}^{-1}}
\newcommand{\ug}[1]{u_{\gamma_{#1}}}

\newcommand{\Supp}{\mathop{\rm Supp}\nolimits}
\newcommand{\FAT}[1]{\mbox{{$\mathbb{#1}$}}}

\newcommand{\zz}{\FAT{Z}}
\newcommand{\cc}{\FAT{C}}
\newcommand{\qq}{\FAT{Q}}

\newtheorem{lem}{\sc Lemma}
\newtheorem{prop}[lem]{\sc Proposition}

\newtheorem{thm}[lem]{\sc Theorem}
\newtheorem{df}[lem]{\sc Definition}
\input{amssym.def}
\input{amssym.tex}
\begin{document}

\title{Projective bases of division algebras and groups of central type II\footnote{This
work was supported in part by the US--Israel
Binational Science Foundation Grant 82334.}}
\author{Michael Natapov \\
\small Department of Mathematics, Technion, Haifa 32000,
Israel\footnote{Current address: Department of Mathematics, Indiana University, Bloomington 47405.} \\
\small e-mail: mnatapov@indiana.edu}
\date{}
\maketitle
\begin{abstract}
Let $G$ be a finite group and let $k$ be a field. We say that $G$
is a projective basis of a $k$--algebra $A$ if it is isomorphic to
a twisted group algebra $k^\alpha G$ for some $\alpha \in
H^2(G,k^\times)$, where the action of $G$ on $k^\times$ is
trivial. In a preceding paper by Aljadeff, Haile and the author it
was shown that if a group $G$ is a projective basis in a
$k$--central division algebra then $G$ is nilpotent and every
Sylow $p$--subgroup of $G$ is on the short list of $p$--groups,
denoted by $\Lambda$. In this paper we complete the classification
of projective bases of division algebras by showing that every
group on that list is a projective basis for a suitable division
algebra.

We also consider the question of uniqueness of a projective basis
of a $k$--central division algebra. We show that basically all
groups on the list $\Lambda$ but one satisfy certain rigidity
property.
\end{abstract}

{\bf 1. Introduction.}

Let $k$ be a field. Let $A$ be a $k$--central simple algebra. A
basis $\{a_1, a_2, \ldots a_n \}$ of $A$ is called a
\emph{projective basis} if all $a_i$-s are invertible and for
every pair $i,j$ there is an $m$ such that $a_ia_j =
\lambda_{ij}a_m$ for some $\lambda_{ij} \in k^\times$. It is not
difficult to see that $A$ has a projective basis if and only if it
is isomorphic to a twisted group algebra $k^\alpha G$ for some
finite group $G$ and some $\alpha \in H^2(G,k^\times)$, where the
action of $G$ on $k^\times$ is trivial. The most important
examples of algebras with projective bases are the symbol
algebras. Recall that a $k$--central simple algebra $B$ of
dimension $n^2$ is a symbol algebra if $B$ is generated by two
elements $x$ and $y$ with relations $x^n \in k^\times$, $y^n \in
k^\times$, $xy = \xi_n yx$ ($\xi_n$ is a primitive $n$--root of
unity contained in $k$). It is easy to see that $B$ is isomorphic
to $k^\alpha (Z_n \times Z_n)$ for a suitable $\alpha \in H^2(Z_n
\times Z_n, k^\times)$, where $Z_n$ denotes  the cyclic group of
order $n$. In fact, if $G$ is abelian and $k^\alpha G$ is a
$k$--central simple algebra then $G$ is of \emph{symmetric type}
(i.e. $\cong H \times H$ for some abelian group $H$) and $k^\alpha
G$ is isomorphic to a tensor product of symbol algebras (see e.g.
\cite[Theorem 1.1]{AS}).

Central simple algebras with projective bases appear in the theory
of $G$--graded algebras. Recall that an (associative) algebra $A$
over a field $k$ is \emph{graded} by a group $G$ if $A$ decomposes
into the direct sum of $k$--vector subspaces $A = \oplus_{g \in G}
A_g$ such that $A_g A_h \subseteq A_{gh}$ for any $g,h \in G$. A
$G$--grading on $A$ is called \emph{fine} if $dim_k(A_g) \leq 1$
for all $g \in G$ (see \cite{BSZ} for more details). Clearly, if
$A$ is isomorphic to a twisted group algebra $k^\alpha G$ then it
is endowed with a fine $G$--grading over $k$. Conversely, it is
shown in \cite[Theorem 1]{AHN2}, that the support $\Supp A = \{g
\in G : A_g \neq 0 \}$ of a fine $G$--grading on a $k$--central
simple algebra $A$ is a projective basis of $A$.

Groups $G$ which are projective bases of central simple algebras
are of special interest in the representation theory of finite
groups. Recall that the dimension of an irreducible representation
of a finite group $\Gamma$ is not greater than the square root of
$[\Gamma:Z(\Gamma)]$, where $Z(\Gamma)$ denotes the center of
$\Gamma$. By definition, the group $\Gamma$ is of \emph{central
type} if it admits an irreducible representation of the maximal
possible dimension $\sqrt{[\Gamma:Z(\Gamma)]}$. A remarkable
result of DeMeyer and Janusz establishes that $\Gamma$ is of
central type if and only if every Sylow $p$--subgroup $S_p$ of
$\Gamma$ is of central type and $Z(S_p)=Z(\Gamma)\bigcap S_p$
(\cite[Theorem 2]{DMJ}). Isaacs and Howlett proved, using the
classification of finite simple groups, that if $\Gamma$ is of
central type then it is solvable (\cite [Theorem 7.3]{HI}).

If $\Gamma$ is a group of central type and $G = \Gamma/Z(\Gamma)$,
then the irreducible representation of $\Gamma$ of dimension
$\sqrt{[\Gamma:Z(\Gamma)]}$ induces a projective representation of
$G$ (of the same dimension), and so there exists a cohomology
class $\alpha \in H^2(G, \cc^\times)$, such that $\cc^\alpha G
\cong M_n(\cc)$. By abuse of language, we will call such $G$ a
group of central type as well. In fact, it is easy to see that a
group $G$ is of central type in this new sense if and only if $G
\cong \Gamma/Z(\Gamma)$, where $\Gamma$ is some group of central
type in the classical sense. Note, that the result of Isaacs and
Howlett holds for a group of central type in the new sense as
well. Also, $G$ is of central type if and only if every Sylow
$p$--subgroup of $G$ is of central type (\cite[Corollary 4]{DMJ}).
In this paper we will use the notion of a group of  central type
only in the new sense.

In \cite{AH} Aljadeff and Haile analyzed division algebras which
contain projective bases $G$. In particular they obtained two
necessary conditions on the group $G$.

\begin{thm} \label{nilpotent} If $k^\alpha G$ is a division algebra with
center $k$ then $G$ is nilpotent and its commutator subgroup is
cyclic \rm{(\cite[Theorems 1 and 2]{AH})}.
\end{thm}

It follows by the nilpotency condition that a $k$--central
division algebra $k^\alpha G$ is isomorphic to
$k^{\alpha_1}P_1\otimes k^{\alpha_2}P_2\otimes\ldots \otimes
k^{\alpha_m}P_m$, where $P_1,P_2,\dots, P_m$ are the Sylow
$p$--subgroups of $G$ and $\alpha_i$ is the restriction of
$\alpha$ to $P_i$. Conversely, if $P_1,P_2,\ldots, P_m$ are
$p$--groups (for $m$ different primes) and $k^{\alpha_i}P_i$ is a
$k$--central division algebra for all $i$, then
$k^{\alpha_1}P_1\otimes k^{\alpha_2}P_2\otimes\ldots \otimes
k^{\alpha_m}P_m$ is a division algebra with projective basis $G
\cong P_1 \times \ldots \times P_m$. This reduces the analysis of
such algebras to the case where $G$ is a $p$--group.

In \cite [Corollary 3]{AHN} there is a (short) list $\Lambda$ of
$p$--groups containing all $p$--groups which possibly are
projective bases of division algebras. The list $\Lambda$ consists
of three families of groups $G$:

\begin{enumerate}
   \item $G$
   is abelian of symmetric type, that is $G \cong \prod(Z_{p^{n_i}} \times Z_{p^{n_i}})$,
   \item $G \cong G_1 \times G_2$ where $G_1 = Z_{p^n} \rtimes Z_{p^n} =
   \langle\pi, \sigma \mid \sigma^{p^n}=\pi^{p^n}=1 \; {\rm and} \;
   \sigma\pi\sigma^{-1} = \pi^{p^s+1}\rangle$ where $1 \leq s < n$ and
   $1 \neq s$ if $p=2$,
   and $G_2$ is an abelian group of symmetric type of exponent
   $\leq p^s$,
   \item $G \cong G_1 \times G_2$ where \\
   $G_1 = Z_{2^{n+1}} \rtimes (Z_{2^n} \times Z_2) =
               \bigg \langle \pi,\sigma,\tau \, \bigg|  \,
               \begin{array}{l}
                 \pi^{2^{n+1}}=\sigma^{2^n}=\tau^2=1, \sigma\tau = \tau\sigma,\\
                 \sigma\pi\sigma^{-1}=\pi^3,
                 \tau\pi\tau^{-1}=\pi^{-1}
               \end{array}\bigg\rangle$ and
   $G_2$ is an abelian group of symmetric type of exponent $\leq 2$.
\end{enumerate}

For the reader convenience we record \cite [Corollary 3]{AHN} in
the following theorem:

\begin{thm} If a $p$--group $G$ is a projective basis of a
division algebra then $G$ is in $\Lambda$.
\end{thm}

The main purpose of this paper is to complete the classification
of projective bases of division algebras, begun in \cite{AHN}, by
showing that every group on the list $\Lambda$ is a projective
basis for a suitable division algebra over an appropriate field.
Thus, combining this with Theorems 1 and 2 we have the following
result:

\begin{thm} \label{classification} Let $G$ be a finite group.
Then there exist a field $k$ and a cohomology class $\alpha \in
H^{2}(G,k^\times)$ such that the twisted group algebra $k^{\alpha
}G$ is a $k$--central division algebra if and only if $G$ is
nilpotent and all Sylow $p$--subgroups of $G$ are in $\Lambda$.
\end{thm}

Now, combining Theorem \ref{classification} with \cite[Theorem
1]{AHN2} we obtain a complete classification of the groups which
support fine gradings on finite dimensional division algebras over
their centers:

\begin{thm} Let $G$ be a finite group.
Then there exist a field $k$ and a $k$--central division algebra
$D$ with a fine grading such that $\Supp D = G$ if and only if $G$
is nilpotent and all Sylow $p$--subgroups of $G$ are in $\Lambda$.
\end{thm}

Next we consider the question of uniqueness of a projective basis
of a $k$--central division algebra.

\noindent {\sc Question.} (Strong rigidity).  \emph{Let $k^{\alpha
}G$ and $k^{\beta }H$ be isomorphic $k$--central division
algebras. Is necessarily $G \cong H$?}

The answer is negative in general. One can build a division
algebra which has two non isomorphic abelian projective bases, see
e.g. the construction in \cite{Set}. Moreover, it is shown in
\cite[proof of Theorem 13]{AHN} that any $k$--central division
algebra of the form $k^\alpha (Z_4 \rtimes (Z_2 \times Z_2))$ is
isomorphic to $k^\beta (Z_2 \times Z_2 \times Z_2 \times Z_2)$ for
a suitable $\beta \in H^2( (Z_2)^{\times 4}, k^\times)$. The
second objective of this paper is to show that the group $Z_4
\rtimes (Z_2 \times Z_2)$ is basically the only group on the list
$\Lambda$ which does not satisfy the following weak version of
rigidity.

\begin{df} We say that a group $G$ satisfies weak rigidity if there exist a
field $k$ and a cohomology class $\alpha \in H^{2}(G,k^\times)$
such that $k^{\alpha }G$ is a $k$--central simple algebra and if
$k^{\alpha }G \cong k^{\beta }H$ for some $H$ and $\beta \in
H^2(H, k^\times)$ then $H \cong G$.
\end{df}

 Our result is given in the following theorem:

\begin{thm} \label{rigidity} If a group $G \in \Lambda$ has no direct factor
isomorphic to ${Z_4 \rtimes (Z_2 \times Z_2)}$, then $G$ satisfies
weak rigidity.
\end{thm}

\bigbreak
\par

{\bf 2. Realization.}

In this section we prove Theorem \ref{classification}. Of course,
we have to show only one direction (namely, the part ``if" of the
theorem). In Case (I) below we exhibit the construction of a
division algebra with projective basis $G = Z_{p^n} \rtimes
Z_{p^n}$ (cf. \cite{AH}, p. 192). Next, in Case (II) we realize
the group $Z_{2^{n+1}} \rtimes (Z_{2^n} \times Z_2)$ as a
projective basis of division algebra, and then, in Case (III), we
show how to realize arbitrary $p$--groups on the list $\Lambda$.
All the realizations are similar and done over the field of
iterated Laurent series $k = K((t_1))((t_2))\ldots ((t_N))$ where
the field $K$ and $N$ depend on $G$.

\textbf{(I)} Let $G=\langle \pi,\sigma \mid
\sigma^{p^n}=\pi^{p^n}=~1,
     \sigma \pi \sigma^{-1} = \pi^{p^s+1} \rangle$, $s \leq n$ and
$s\neq1$ if $p=2$. (It is an abelian group when $s = n$).

Let $K$ be a field of characteristic zero that contains a
primitive $p^s$--root of unity $\xi$ and does not contain
primitive $p^{s+1}$ roots of unity. For any $c \in K^\times$, let
$L = K(\up)/K$ be a cyclic Galois extension defined by $\up^{p^n}
= c^{p^{n-s}} \xi$ with the Galois group $Gal(L/K) \cong Z_{p^n}$.
Since $c^{-1}\up^{p^s}$ is a primitive $p^n$--root of unity, a
generator $\sigma$ of the Galois group of $L$ can be chosen such
that $\sigma(\up) = c^{-1}\up^{p^s+1}$.

Let $t$ be an indeterminate and let $k = K((t))$ be the field of
iterated Laurent series over $K$. Consider the field $L((t)) \cong
L \otimes_K K((t))$ which is a cyclic extension of $K((t))$ with
the same Galois group $\langle \sigma \rangle$ as that of $L$.
Consider the cyclic crossed product $D = (L((t))/k, \sigma, t)$,
that is $D = \oplus_{i=0}^{p^n-1} L((t)) \us^i$ as an
$L((t))$--vector space with the multiplication given by
$\us b = \sigma(b) \us$ for any $b \in L((t))$ and $\us^{p^n} =
t$. We claim that the group $G$ is a projective basis of $D$.
Indeed, let $\Gamma$ denote the multiplicative subgroup of
$D^\times$ generated by $\up$ and $\us$. Then $D = k(\Gamma)$,
that is $D$ is generated by $\Gamma$ as a $k$--vector space. It is
easy to see that $k^\times \Gamma /k^\times \cong G$ and $D \cong
k^\alpha G$ where $\alpha$ corresponds to the central extension
$1\rightarrow k^\times \rightarrow k^\times \Gamma \rightarrow G
\rightarrow 1$. Observe that in case that $G$ is abelian, the
algebra constructed above is isomorphic to the symbol algebra $(c
\xi, t)_{p^n}$.

Now, it is well known that $D$ is a division algebra. One way to
show this is to view $D$ as a ring of twisted Laurent series over
the field $L$ in the variable $\us$. Namely, $D =  L((\us; \sigma))
= \left\{ \sum_{i\geq k}a_i\us^i \mid k \in \zz, a_i \in L \right\}$
with the multiplication on $L((\us; \sigma))$ given by $\us b =
\sigma(b) \us$ for any $b \in L$, where $\sigma$ is the automorphism
of $L$ defined above. This proves that $D$ is a division algebra
(see \cite[Example 1.8]{L}). We use this argument in Cases (II) and
(III) below.

 \textbf{(II)} $G = Z_{2^{n+1}} \rtimes (Z_{2^n} \times Z_2)=
                    \bigg \langle \pi,\sigma,\tau \, \bigg|
                    \begin{array}{l}
                     \pi^{2^{n+1}}=\sigma^{2^n}=\tau^{2}=1, \sigma \tau = \tau
                     \sigma, \\
                     \sigma\pi\sigma^{-1}={\pi^3}, \tau\pi\tau^{-1}=\pi^{-1}
                    \end{array}\bigg\rangle$.

Let $K$ be a field of characteristic zero that does not contain
$\sqrt{-1}, \sqrt{2}$ and $\sqrt{-2}$.  For any $c \in K^\times$
such that $c^{2^n} \notin 4K^4$, we let $L = K(\up)/K$ be a Galois
extension defined by $\up^{2^{n+1}} = -c^{2^n}$. The Galois action
of $Gal(L/K) \cong Z_{2^n} \times Z_2 = \langle \sigma, \tau
\rangle$ on $L$ is given by
\[
 \sigma(\up) =  c^{-1}\up^3 \;\;\; {\rm and} \;\;\; \tau(\up) =
 c\Up.
\]

Let $D_1 = L((\us; \sigma))$ be a ring of twisted Laurent series
over $L$ in a variable $\us$. As above, it is a division algebra.
Next, let $D = D_1((\ut; \tau))$ be a ring of twisted Laurent series
over the algebra $D_1$ in a variable $\ut$, where the automorphism
$\tau$ of $D_1$ extends the action of $\tau$ on $L$ and the action
on $\us$ is trivial. Since $D_1$ is a division algebra, $D$ is a
division algebra as well.

It is easy to see that the center $k$ of $D$ is generated by the
field $K = L^{\langle \sigma, \tau \rangle}$ and the elements $s =
\us^{2^n}$ and $t = \ut^2$, namely $k = K((s))((t))$. Moreover, the
field $L((s))((t))$ which is a Galois extension of $k$ with the
Galois group $\langle \sigma, \tau \rangle$, is a maximal subfield
of $D$. That is $D$ is isomorphic to the crossed product
$(L((s))((t))/k, \langle \sigma, \tau \rangle, f)$.  The elements
$\us$ and $\ut$ represent $\sigma$ and $\tau$ in $D$ and the
2-cocycle $f$ is given by
\[
 \us^{2^n} = s, \;\;\; \ut^2 = t \;\;\; {\rm and} \;\;\;  (\us, \ut) = 1,
\]
where $(\us, \ut)$ denotes the commutator of $\us$ and $\ut$.
Finally, arguing as in the previous case we see that $D$ is
isomorphic to a twisted group algebra $k^\alpha G$ for an
appropriate class $\alpha \in H^2(G, k^\times)$.

\textbf{(III)} We complete the realization of $p$--groups as
follows. Let $G$ be a group on the list $\Lambda$. Write $G = G_0
\times Z_{p^r} \times Z_{p^r}$. We assume, by induction, that the
subgroup $G_0$ is realizable as a projective basis of a division
algebra, namely, there exist a field $K$ and a cohomology class
$\beta \in H^2(G,K^\times)$, such that $D_0 = K^\beta G_0$ is a
division algebra. We may assume also that $K$ contains a primitive
$p^r$--root of unity. Let $k = K((s))((t))$ where $s,t$ are
indeterminates, and consider the $k$--algebra $D = k^\beta G_0
\otimes_k (s,t)_{p^r}$, where $k^\beta G_0 \cong K^\beta G_0
\otimes_K k$ and $(s,t)_{p^r}$ is a symbol algebra. Clearly, $D
\cong k^\alpha G$ for some $\alpha \in H^2(G,k^\times)$, such that
$res^G_{G_0} (\alpha) = \beta$, that is $G$ is a projective basis
of $D$.

We now show that $D$ is a division algebra. Let $x$ and $y$ be
standard generators of the symbol $(s,t)_{p^r}$, that is $x^n=s$,
$y^n=t$ and $xy = \zeta yx$ ($\zeta$ is a primitive $p^r$--root of
unity). Let $D_1 = D_0((x))$ be a ring of Laurent series in the
variable $x$ over $D_0 = K^\beta G_0$. Since $D_0$ is a division
algebra it follows that $D_1$ is a division algebra as well. Now,
let $D_1((y; \tau))$ be a twisted Laurent series ring over $D_1$ in
the variable $y$, where the automorphism $\tau$ of $D_1$ is trivial
on $D_0$ and $\tau(x) = \zeta^{-1} x$. Clearly, $D \cong D_1((y;
\tau))$ and hence $D$ is a division algebra.

\textbf{(IV)} Now, let $G$ be a nilpotent group. Write $G$ as a
direct product $G = P_1 \times \ldots \times P_m$ of its Sylow
$p_i$--subgroups. Suppose $P_i \in \Lambda$ for all $1\leq i\leq
m$. For every $i$, we can construct as above a field $k_i$ and a
cohomology class $\alpha_i \in H^2(P_i, k_i^\times)$ such that
$k_i^{\alpha_i} P_i$ is a division algebra. Moreover, we can
choose the field $k_i$ to be $\qq(\xi_i)((t_1))\ldots
((t_{N_i}))$, where $\xi_i$ is a $p_i^{s_i}$--primitive root of
unity, for a suitable number of indeterminates $N_i$. Let $\xi =
\prod_i \xi_i$ be a root of unity of order $\prod_i p_i^{s_i}$ and
let $K_i = \qq(\xi)((t_1))\ldots ((t_{N_i}))$. Observe that since
$K_i$ does not contain $p_i^{s_i + 1}$--primitive roots of unity,
precisely the same construction of $\alpha_i \in H^2(P_i,
K_i^\times)$ as in {(I - III)} gives a division algebra
$K_i^{\alpha_i} P_i$. Consider $K_1,\ldots,K_m$ as subfields of $k
= \qq(\xi)((t_1))\ldots ((t_{N}))$ where $N = \max_i(N_i)$. For
all $1\leq i\leq m$, let $D_i = K_i^{\alpha_i} P_i \otimes k$. By
\cite[Corollary 19.6 a]{Pi}, we see that $D_i$ is a division
algebra. Finally, $D = D_1 \otimes_k \ldots \otimes_k D_m$ is a
division algebra, since all $D_i$ have relatively prime degrees,
and $G$ is a projective basis of $D$.

This completes the proof of Theorem \ref{classification}.\qed

We close this section by pointing out that in all of our
constructions we may replace Laurent series by rational functions.
Indeed, given a group $G$ as in (IV), we may follow the above
construction but now over an appropriate field of rational
functions of the form $\qq(\xi)(t_1, \ldots, t_{N})$, to obtain a
central simple algebra $A$. This algebra restricted to the Laurent
series field $\qq(\xi)((t_1))\ldots ((t_{N}))$ is a division
algebra and therefore $A$ is a division algebra as well.

\bigbreak
\par
{\bf 3. Rigidity.}

In this section we prove Theorem \ref{rigidity}.

We first prove the theorem for abelian $p$--groups. Let $G$ be an
abelian group of symmetric type, that is $G = \prod_{k=1}^{\ell}
Z_{p^{n_k}}\times Z_{p^{n_k}}$. We construct a division algebra
$D$ such that any projective basis of $D$ is isomorphic to $G$.
Let $F = \cc((t_1))\ldots ((t_N))$, $N \geq 2\ell$, denote the
$N$--fold iterated Laurent series field over $\cc$ (the Amitsur
field). Consider the set of symbol algebras $\{
(t_{2k-1},t_{2k})_{p^{n_k}} \}_{k=1}^\ell$ over the field $F$, and
let $i_k, j_k$ be their standard generators ($i_k$ and $j_k$
satisfy $i_k^{p^{n_k}} = t_{2k-1}$, $j_k^{p^{n_k}} = t_{2k}$ and
$i_k j_k = \xi_{p^{n_k}} j_k i_k$ where $\xi_{p^{n_k}}$ is a
primitive $p^{n_k}$--root of unity). Let
\begin{equation} \label{valued}
 D =\bigotimes_{k=1}^{\ell} (t_{2k-1},t_{2k})_{p^{n_k}}.
\end{equation}
Clearly $D$ is isomorphic to a twisted group algebra $F^\alpha G$
for an appropriate class $\alpha \in H^2(G, F^\times)$. Moreover,
$D$ is a division algebra by \cite[Example 3.6 (a)]{TW}.

\begin{prop} \label{abelian} Let $G$, $F$ and $D \cong F^\alpha G$
be as above. Let $H$ be a group and $\beta \in H^2(H, F^\times)$.
If $D \cong F^\beta H$ then $G \cong H$.
\end{prop}

In order to prove the Proposition we view $D$ as a valued tame and
totally ramified (TTR) division algebra over $F$.

Let us recall some definitions and notation related to valuations
on division algebras (cf. \cite{TW}). Let $v$ be a valuation on an
$F$--central division algebra $D$ with values in a totally ordered
abelian group $\Gamma$. We let $\Gamma_D = v(D^\times)$ and
$\Gamma_F = v(F^\times)$ be the value group of $v$ on $D$ and $F$,
respectively. The algebra $D$ is called \emph{tame and totally
ramified} over $F$ with respect to $v$ if $|\Gamma_D: \Gamma_F| =
[D:F]$ and $char(\overline{F}) \nmid [D:F]$, where $\overline{F}$
is the residue class field of $F$.

We now define a valuation on the Amitsur field $F =
\cc((t_1))\ldots ((t_N))$ and its extension to the division
algebra $D$ defined in (\ref{valued}). Consider the group $\zz^N$
with the right-to-left lexicographic order. There is a valuation
$v$ on $F$ with values in $\zz^N$:
\[
v\left( \sum_{i_1}\cdots \sum_{i_N} c_{i_1\ldots i_N}
t_1^{i_1}\cdots t_N^{i_N}\right) = min\{ (i_1,\ldots, i_N) \; | \;
c_{i_1\ldots i_N} \neq 0 \}.
\]
The valuation $v$ is called the \emph{standard} valuation on $F$.
Its value group is $\Gamma_F = \zz^N$ and its residue field is
$\overline{F} = \cc$.

The division algebra $D$ defined in (\ref{valued}) has a valuation
$v: D^\times \rightarrow \qq^N$ which extends the standard
valuation $v$ on $F$:
\begin{equation}  \label{valuation}
  v(i_k) = \frac{1}{p^{n_k}}v(t_{2k-1}) =
(0,\ldots,0,\frac{1}{p^{n_k}},0,\ldots,0),
\end{equation}
\[  v(j_k) =
\frac{1}{p^{n_k}}v(t_{2k})=(0,\ldots,0,\frac{1}{p^{n_k}},0,\ldots,0),
\]
(with nonzero entries in the $2k-1$ and $2k$ positions
respectively). With respect to the valuation $v$ we have $\Gamma_D
= \langle v(i_1),v(j_1),\ldots, v(j_\ell) \rangle + \Gamma_F$, and
so $\Gamma_D/ \Gamma_F$ (the \emph{relative value group} of $D$
with respect to $v$) is isomorphic to $\prod_{k=1}^{\ell}
Z_{p^{n_k}}\times Z_{p^{n_k}} (\cong G)$. Therefore, the division
algebra $D$ is TTR over $F$.

Next, we recall the notion of \emph{armature} (\cite{TW}) which is
basically the same as the notion of an abelian projective basis:

\begin{df} Let $A$ be a finite-dimensional $F$--algebra.
Let $\mathcal{A}$ be a (finite) subgroup of $A^\times/F^\times$
and $a_1, a_2, \dots a_n$ be a representatives of the elements
$\mathcal{A}$ in $A$. We say $\mathcal{A}$ is an armature of $A$
if $\mathcal{A}$ is abelian and $\{ a_1, a_2, \dots a_n \}$ is an
$F$--base of $A$.
\end{df}

Clearly, the group generated by $\{ i_kF^\times/F^\times,
j_kF^\times/F^\times \}_{k=1}^\ell$ in $D^\times/F^\times$ is an
armature of $D$. The following result (due to Tignol and
Wadsworth, \cite[Proposition 3.3]{TW}) establishes that the
armature of the algebra $D$ is uniquely determined by its relative
value group.

\begin{prop} \label{armature}Let $(D,v)$ be a valued division algebra
with $D$ tame and totally ramified over its center $F$. If
${\mathcal{A}}$ is an armature of $D$ as an $F$--algebra then the
map $\overline{v}: {\mathcal{A}} \rightarrow \Gamma_D/\Gamma_F$
induced by $v$ is an isomorphism.
\end{prop}

Now, we can prove Proposition \ref{abelian}.

\begin{proof} Let $H$ be an abelian group and suppose there exists
a cohomology class $\beta \in H^2(H, F^\times)$ such that $F^\beta
H \cong D$. Note that $\mathcal{B} = \langle \us F^\times/F^\times
\; | \; \sigma \in H \rangle \cong H$ is an armature of $D$. By
Proposition \ref{armature}, $\mathcal{B}$ is isomorphic to the
relative value group $\Gamma_D/\Gamma_F$ with respect to the
valuation $v$ defined in (\ref{valuation}). Since
$\Gamma_D/\Gamma_F \cong G$, we get $G \cong H$.

A nonabelian group cannot form a projective basis of a division
algebra over the Amitsur field $F$, because $F$ contains all roots
of unity (see \cite[Sec. 2]{AH}). Hence the algebra $D$ has no
nonabelian projective basis and the proposition follows.
\end{proof}

It remains to prove Theorem \ref{rigidity} for nonabelian groups.

\medskip Case I. $G = (Z_{p^n} \rtimes Z_{p^n}) \times Z_{p^{r_2}} \times
Z_{p^{r_2}} \times \ldots \times Z_{p^{r_\ell}} \times
Z_{p^{r_\ell}}$ with a set of generators $\pi, \sigma, \gamma_3,
\ldots , \gamma_{2\ell}$. Assume that $G' = \langle \pi^{p^s}
\rangle$ ($s \geq 1$, or $s \geq 2$ when $p = 2$) and (therefore)
$r_k \leq s$ for all $2\leq k \leq \ell$.

For $N = 2\ell$ define $K = \qq (\xi) ((t_1))\ldots ((t_N))$ where
$\xi$ is a primitive $p^s$--root of unity. It was shown in the
previous section (see (I), (III)) that there is a class $\alpha
\in H^2(G, K^\times)$ such that $D \cong K^\alpha G$ is a division
algebra. Namely, we let $D$ be a tensor product of the form $D =
D_1 \otimes D_2 \otimes \ldots \otimes D_\ell$, where $D_1$ is a
cyclic algebra generated by the elements $\up$ and $\us$ subject
to the following relations:
\[
 \up^{p^n} = t_1^{p^{n-s}}\xi, \;\;\;
 \us^{p^n} = t_2       \;\;\; {\rm and} \;\;\;
 (\us, \up) = t_1^{-1} \up^{p^s},
\]
and for all $2\leq k \leq \ell$, $D_k$ is the symbol algebra
$(t_{2k-1},t_{2k})_{p^{r_k}}$.

We \underline{claim} that the algebra $D$ is of exponent $p^n$.
Indeed, the algebra $D_1$ is isomorphic to a cyclic algebra of the
form $\big(\qq (\xi) ((t_1))(\up)((t_2)), \sigma, t_2 \big)
\otimes_{\qq (\xi) ((t_1))((t_2))} K$ and it is of exponent $p^n$ by
\cite[Corollary 19.6 b]{Pi} and \cite[Corollary 19.6 a]{Pi}.
Furthermore, the symbol $D_k$ is of exponent $p^{r_k} \leq p^s <
p^n$ for all $k$, and the claim follows.

Suppose that $D$ is isomorphic to $K^\beta H$ for some $H$ and
$\beta \in H^2(H, K^\times)$. Observe that $H$ is not abelian, for
otherwise, by \cite[Theorem 1.1]{AS}, $K^\beta H$ is a tensor
product of symbol algebras, and hence $exp(K^\beta H) \leq p^s$ --
the number of $p$--power roots of unity in the field $K$, a
contradiction. Therefore, by \cite[Theorem 1]{AHN}, $H$ is of the
form $(Z_{p^m} \rtimes Z_{p^m}) \times B$ where generators $x$ and
$y$ of the semidirect product $Z_{p^m} \rtimes Z_{p^m}$ satisfy
$x^{p^m} = y^{p^m} = 1$ and $yxy^{-1} = x^{p^s + 1}$ and $B =
Z_{p^{f_1}} \times Z_{p^{f_1}} \times \ldots \times
Z_{p^{f_\jmath}} \times Z_{p^{f_\jmath}}$ is abelian of symmetric
type of exponent $\leq p^s$.

Consider the subalgebra $K^\beta (Z_{p^m} \rtimes Z_{p^m})$ of
$K^\beta H$. By the Factorization Lemma in \cite{AH} it can be
factored from $K^\beta H$, that is there exists a 2--cohomology
class ${\widetilde{\beta}}$ on $B \cong H/(Z_{p^m} \rtimes
Z_{p^m})$ such that:
\[
 K^\beta H \cong K^\beta (Z_{p^m} \rtimes Z_{p^m}) \otimes
 K^{\widetilde{\beta}}
 B.
\]
Since $B$ is abelian, by \cite[Theorem 1.1]{AS},
$K^{\widetilde{\beta}} B$ is a product of symbol algebras of the
form:
\[
K^{\widetilde{\beta}} B = \bigotimes_{k=1}^{\jmath}
(a_{2k-1},a_{2k})_{p^{f_k}}.
\]
In particular, it follows that the algebra $K^\beta H$ is of
exponent at most $p^m$.

We \underline{claim} that $m = n$. First, if $m < n$ then
$exp(K^\beta H) \leq p^m < p^n = exp(D)$, a contradiction. To see
that $m \leq n$, we restrict $D = D_1 \otimes \ldots \otimes
D_\ell$ to the Amitsur field $F = \cc((t_1))\ldots ((t_N)) \cong K
\otimes_{\qq(\xi)} \cc$.

Consider the subfield $E = K(z)$ of $D_1$ where $z =
u_{\pi}^{p^s}/t_1$. By \cite[Proposition 7.2.2]{R} $D_1 \otimes E$
is Brauer equivalent to the centralizer $C_{D_1}(E)$ of $E$ in
$D_1$. It is easy to see $C_{D_1}(E) = K(u_\pi, \us^{p^{n-s}})$,
and it is isomorphic to the symbol algebra $(t_1 z, t_2)_{p^s}$
over the field $E$. Since $z = u_{\pi}^{p^s}/t_1$ is a primitive
$p^{n}$--root of unity, it follows that $D_1 \otimes F \cong D_1
\otimes_K K(\zeta) \otimes_{K(\zeta)} F \sim (t_1 \zeta,
t_2)_{p^s}$ where $\zeta$ is a primitive $p^{n}$--root of unity in
$\cc$. Since the symbol algebra $(t_1 \zeta, t_2)_{p^s}$ is Brauer
equivalent to $(t_1,t_2)_{p^s}$ over $F$ (\cite[Proposition
7.1.17]{R}) we have:
\begin{equation} \label{restricted}
D \otimes_K F  \sim  (t_1,t_2)_{p^s} \otimes_F (t_3,t_4)_{p^{r_2}}
\otimes_F \ldots \otimes_F (t_{2\ell-1},t_{2\ell})_{p^{r_\ell}}.
\end{equation}
Since the latter is a TTR division algebra we have that $Ind(D
\otimes F) = \frac{Ind(D)}{p^{n-s}}$. Now consider the
multiplicative subgroup $\mathcal{H}$ of $(K^\beta H)^\times$
generated by representatives of $H$ in $K^\beta H$. Observe that
$\mathcal{H}$ is center by finite, so by a theorem of Schur
\cite[Chapter 2, Theorem 9.8]{Su} its commutator subgroup
$\mathcal{H}'$ is finite. It is easy to see that $K^\times
\mathcal{H}'/K^\times = H'$. Since the commutator subgroup $H'$ of
$H$ is of order $p^{m-s}$ it follows that $K^\beta H$ contains a
cyclotomic field extension $K^\beta H'/K$ of degree $p^{m-s}$ and
hence $Ind(K^\beta H \otimes F) \leq \frac{Ind(K^\beta
H)}{p^{m-s}}$. Thus we have $m \leq n$ and the claim follows.

Now, write $H = (Z_{p^n} \rtimes Z_{p^n}) \times Z_{p^{f_1}}
\times Z_{p^{f_1}} \times \ldots \times Z_{p^{f_\jmath}} \times
Z_{p^{f_\jmath}}$, and let $x$ and $y$ be generators of the
semidirect product $Z_{p^n} \rtimes Z_{p^n}$. Let $u_x,
u_{x^{p^s}}$ be representatives of $x$ and $x^{p^s}$ in $K^\beta
H$. Since the field $K^\beta H' = K(u_{x^{p^s}})$ is a cyclotomic
extension of $K$, we may assume that $u_{x^{p^s}}^{p^{n-s}} =
\xi$. There is an element $a \in K^\times$ such that $u_x^{p^s} =
a u_{x^{p^s}}$. It follows that $u_x^{p^n} = a^{p^{n-s}}\xi$.
Since $\langle x \rangle$ is a normal subgroup of $H$, by
\cite[Lemma A]{AH} we have that $K(u_x)/K$ is a Galois field
extension which is cyclic of order $p^n$. Moreover, conjugation by
representatives $u_h, h \in H$ of $K^\beta H$ induces a surjective
homomorphism $H/\langle x \rangle \rightarrow Gal(K(u_x)/K)$. It
follows that conjugation by a representative $u_y$ of $y$ induces
a Galois action on $K(u_x)$, and we may assume (choosing a new
generator $y$ if necessary) that $u_y u_x u_y^{-1} = a^{-1}
u_x^{p^s + 1}$. Also, there is an element $b \in K^\times$ such
that $u_y^{p^m} = b$. As in the claim above we have that $K^\beta
(Z_{p^n} \rtimes Z_{p^n}) \otimes F$ (where $F$ is the Amitsur
field defined above) is similar to the symbol algebra
$(a,b)_{p^s}$. Hence, by an index argument we have that
$(a,b)_{p^s} \otimes \bigotimes_{k=1}^{\jmath}
(a_{2k-1},a_{2k})_{p^{f_k}}$ is a division algebra, and,
furthermore, it is isomorphic to the algebra obtained in
(\ref{restricted}). Then, applying Proposition \ref{abelian}, we
get $Z_{p^s} \times Z_{p^s} \times Z_{p^{r_2}} \times
Z_{p^{r_2}} \times \ldots \times Z_{p^{r_\ell}} \times
Z_{p^{r_\ell}} \cong Z_{p^s} \times Z_{p^s} \times B$ and hence $G
\cong H$ as well.

\medskip Case II. $G = (Z_{2^{n+1}} \rtimes (Z_{2^n} \times Z_2))
\times Z_2 \times Z_2 \times \ldots \times Z_2\times Z_2$ (where
$n > 1$) with a set of generators $\pi, \sigma, \tau, \gamma_4,
\ldots, \gamma_{2\ell +1}$.

We define $K = \qq((t_1))\ldots ((t_N))$, with $N = 2\ell +1$, and
construct a division algebra $D \cong K^\alpha G$ as follows (see
(II) of the previous section):
\[
D = D_1 \otimes (t_4, t_5) \otimes \ldots \otimes (t_{2\ell},
t_{2\ell +1})
\]
where $D_1$ is generated by elements $\up, \us$ and $\ut$
satisfying the following relations:
\[
 \up^{2^{n+1}} = -t_1^{2^n}, \;\;
 \us^{2^n} = t_2,  \;\; \ut^{2} = t_3, \;\; (\us, \up) = t_1^{-1}\up^2,
 \;\; (\ut, \up) = t_1\up^{-2}, \;\; (\us, \ut) = 1,
\]
and for all $2\leq k \leq \ell$, $(t_{2k}, t_{2k +1})$ is a
quaternion algebra with the standard generators $\ug{2k},
\ug{2k+1}$.

We \underline{claim} that the exponent of $D$ is equal $2^n$.
Indeed, by \cite[Theorem 13]{AHN}, the algebra $D_1$ is isomorphic
to a tensor product of two cyclic algebras, namely
\[
D_1 \cong (K(\up)^\tau, \sigma, t_2) \otimes C
\]
where $C$ is a quaternion algebra. Using the arguments of Case I
we get that the exponent of the cyclic algebra $(K(\up)^\tau,
\sigma, t_2)$ is $2^n$. Thus the claim follows.

Suppose that $D \cong K^\beta H$ for some $H$ and $\beta \in
H^2(H, K^\times)$. Arguing as in Case I, we conclude that the
group $H$ is of the form $(Z_{2^{m+1}} \rtimes (Z_{2^m} \times
Z_2)) \times Z_2 \times Z_2 \times \ldots \times Z_2\times Z_2$.

First, we have $m \geq n$, since by \cite[Theorem 13]{AHN} $K^\beta H$
is isomorphic to a tensor product of cyclic algebras of degrees $2^m$
and $2$, and hence $K^\beta H$ is of exponent at most $2^m$.
Next, we prove that $m \leq n$. Consider the field $E = K(z)$
where $z = u_{\pi}^2/t_1$ is a primitive $2^{n+1}$--root of unity
contained in $D$. We \underline{claim} that $E$ is the maximal
cyclotomic subfield of $D$. Indeed, the centralizer $C_D(E)$ of
$E$ in $D$ is easily seen to be $C_D(E) = K(u_\pi, \us^{2^{n-1}},
\ug{4}, \ldots, \ug{2\ell +1})$, and it is of the form
$$C_D(E) \cong (t_1 z, t_2) \otimes_E (t_4, t_5) \otimes \ldots
\otimes (t_{2\ell}, t_{2\ell +1}),
$$ where the quaternion algebras above are considered over the
field $E$. Now, arguing as in the previous case we see that $D
\otimes F$ and $D \otimes E \sim C_D(E)$ are of the same index and
the claim follows. On the other hand, $K^\beta H$ contains the
cyclotomic extension $K^\beta H'$ of $K$ and its degree is
$ord(H') = 2^m$. This shows that $m \leq n$. Thus we have $m = n$,
and hence $G \cong H$.

This completes the proof of Theorem \ref{rigidity}.\qed

\bigbreak
\par
{\bf Acknowledgments.}

I would like to thank my Ph.D. supervisor at the Technion Eli
Aljadeff for permanent support during the studies, the research and
the work on preparing this paper. I thank Eric Brussel for valuable
discussion of the rigidity question, Darrell Haile and Jack Sonn for
reading the earlier versions of this paper and for many helpful
remarks. Finally I would like to thank the referee for valuable
comments and, in particular, for directing me to the significantly
easier proof of Theorem \ref{classification}.

\bigskip

\end{document}